\input amstex
\documentstyle{amsppt}
   
\loadeufb
\loadeusb
\loadeufm
\loadeurb
\loadeusm

\magnification =\magstephalf
\refstyle{A}
\NoRunningHeads

\topmatter
\title  Uniformization \endtitle
\author Robert  Treger \endauthor
\address Princeton, NJ 08540 \endaddress
\email roberttreger117{\@}gmail.com \endemail
\keywords  Projective variety, holomorphically convex bounded domain
\endkeywords

\document

\head Uniformization 
\endhead

\head
{\rm Robert  Treger} 
\endhead

The main result of this note is the following uniformization theorem.
\proclaim{Theorem}Let $X$ be a connected nonsingular complex projective variety of dimension $n$
without elliptic curves, with ample canonical bundle, and very large residually finite fundamental
group. Then its universal covering $\tilde X$  is a bounded holomorphically convex domain in $\bold
C^n.$
\endproclaim

The uniformization problem for higher-dimensional varieties was  proposed by
Weierstrass
\cite{W, pp. 95, 232, 304} and Hilbert (22nd problem). The above theorem can be viewed as a
converse to the Poincar\'e ampleness theorem
\cite{K, Theorem 5.22}.

We will  describe the idea of the proof.  We take an appropriate Zariski open subset
$Z\subset X$ whose universal covering is a bounded domain $\tilde Z \subset \bold C^n.$ We
construct a one-to-one holomorphic map from the inverse image of $Z$ on $\tilde X$  into $\tilde
Z.$ Finally, we extend the latter inclusion to a holomorphic inclusion
$\tilde X \hookrightarrow \bold C^n.$ 
\head 1. Preliminaries 
\endhead

\subhead{1.1. Definitions and assumptions}\endsubhead Let ${\bold
P}^r=\roman {Proj }\,{\bold C}[z_0, z_1, \dots, z_r] $ be a projective space.
Unless stated otherwise, we
assume $X\subset {\bold
P}^r$ is a connected nonsingular complex projective variety of dimension $n.$ We denote by
$\tilde X$ its universal covering. The fundamental group
$\pi_1(X)$ is called {\it large} \/ \rm if $\tilde X$ contains no proper holomorphic
subsets of strictly positive dimension or, what is the same according to Koll\'ar, for every normal
cycle
$w: W\rightarrow X$, the group
$\roman {im}[\pi_1(W) \rightarrow \pi_1(X)]$
is infinite \cite {K}.

Unless stated otherwise, we also
assume $X$ has no elliptic curves, its  canonical bundle, denoted by  $\eusb K_X,$   is
ample, and  
$\pi_1(X)$ is residually finite and large.   The definition of {\it very large}\/ $\pi_1(X)$ is given in
(2.2.2).

A holomorphic space is called {\it holomorphically convex} \/  if
for any given infinite discrete subset $D$,  there is a holomorphic function on the space which is
unbounded on $D$. It is well known that a holomorphic space is {\it Stein}\/ if and only if it is
holomorphically convex and contains no proper holomorphic subsets of strictly positive dimension.

\subhead{1.2. Grffiths' theorem}\endsubhead One of the main ingredients in the proof of our
theorem is the well-known theorem of Griffiths to the effect that given a point $Q$ on a
nonsingular quasiprojective variety $X
\subset{\bold P}^r $, we may choose a Zariski open neighborhood of $Q$ whose universal
covering is a bounded domain in $\bold C^n$ \cite {G}.

Further, it follows from \cite {G} that one can choose a finite affine Zariski covering of
$X$, denoted by $\{Z_1, \dots , Z_u \}$, such that the universal covering $\tilde Z_i$ of $Z_i$
is a bounded domain and all
$\pi_1 (Z_i)$ are isomorphic residually finite groups without torsion and have the same closed,
normal, residually finite subgroups (identified under the isomorphism). This follows from the
Griffiths theorem since we have a great choice of the neighborhoods while each $\pi_1(Z_i)$ is a
finitely generated group. Furthermore, we can assume that each $Z_i$ is a complement of a 
hyperplane section of $X$ (by taking a Veronese map if necessary).

\subhead{1.3. Fubini spaces}\endsubhead The space $\bold C^{r+1}$, with coordinates $z_0, z_1,
\dots, z_r$, will be considered as a standard Euclidean space with a norm $\root\of {|z_0|^2+
\cdots +|z_r|^2}$. Given a Hilbert space
$H$, we denote by $H^{\ast}$ its dual and by
${\bold P}(H^{\ast})$ the corresponding projective space with the standard
Fubini-Study metric \cite{C, Chap.\, 4}.

For every $i$, we will fix a countable
cofinal set $\{Z_{i,\gamma}
\}_{\gamma\in \bold N}$ of finite Galois coverings of $Z_i$ with the Galois groups independent of\/
$i$ 
$(1\leq i\leq u)$. By the Riemann existence theorem, each covering $Z_{i,\gamma}$ is naturally
compactified to a normal projective variety
$M_{i,\gamma}$ and a finite morphism $\phi_{i,\gamma}: M_{i,\gamma} \rightarrow X$.  
We  assume  the  embedding 
$$
Z_i\subset X \hookrightarrow \bold P^{r} 
$$ 
is given by a very  ample sheaf $\Cal O_X(1)= \eusb K_{X}^{\otimes m}$ ($m \gg 0$). 
The {\it metric}\/ and {\it measure} on $X$ and $Z_i$ are
determined uniquely by this embedding. Here we consider $\bold P^{r}$ with the standard
Fubini-Study metric given by the potential
$$
\log(1 + \sum_{t =1}^r |\zeta_t(x)|^2), \qquad x\in X,
$$
in the affine coordinate system $\zeta_1= z_1\slash {z_0},\; \dots\; , \zeta_r=z_r
\slash {z_0}.$
 The volume form on $X$ can be written locally 
$
dv = \rho \cdot (\root\of {-1})^ndx_1\wedge d{\overline x}_1 \dots \wedge dx_n
\wedge d{\overline x}_n 
,$ 
where $\rho=\rho_X(x)$ is a positive function and $x_1,\dots, x_n$ are local coordinates on $X.$
On finite Galois coverings, it induces volume forms  and  we employ the same notation, $\rho$ and
 $dv.$ 

\subhead{1.4. Positive reproducing kernels and Bergman pseudometrics}\endsubhead Let $U$
denote an arbitrary  complex manifold. Let $B(z,w)$ be a
Hermitian positive definite complex-valued  function  on
$U\times U$ which means:
\roster
 \item"{(i)}" $\overline{B(z, w)}=B(w,z), \qquad B(z,z)\geq 0;$
\item"{(ii)}" $\forall  z_1, \dots, z_N \in U, \quad
 \forall  a_1, \dots, a_N \in \bold C \implies 
\sum_{j,k}^N B(z_k,z_j)a_j \bar {a}_k \geq 0.$
\endroster
 If $B(z, w)$ is, in addition, holomorphic in the first variable then
$B$ is the reproducing kernel of a {\it unique Hilbert space}\/ $H$ of holomorphic functions on $U$ 
(see Aronszajn \cite{A} and the articles by Faraut and Kor\'anyi in \cite{FK, pp.\, 5-14,
pp.\,187-191}). 
The evaluation at a point $Q\in U,$
$
e_Q:  f\mapsto f(Q),
$
is a continuous linear functional on $H$.

Conversely, given a Hilbert space $H$ of holomorphic functions
on
$U$ with all evaluation maps continuous linear functionals then, by the Riesz
representation theorem, for every $w\in U$ there exists a unique function $B_w\in H$
 such that $f(w)= (f, B_w)$ ($\forall f\in H$) and $B(z,w):= B_w(z)$ is the reproducing kernel for
$H$ (which is Hermitian positive definite).

If we assume, {\it in addition}, that $B(z,z) > 0$ for every $z$, then  we can define   
$\log B(z, z)$ and a positive semidefinite  Hermitian form, called the Bergman pseudometric
$$
\qquad \qquad ds^2_U =2 \sum g_{j k} dz_j d\overline z_k, \qquad g_{j
k}={\partial^2\over\partial z_j
\partial \overline z_k  }\log B(z, z) .
$$
 We get a natural map 
$$
\Upsilon : U \longrightarrow \bold {P}(H^*).
$$

In Section 2.2, the function $B(z,w)$ is replace  by a
section of a relevant bundle.

\head  2. Reproducing kernels, metrics, and holomorphic convexity
\endhead
\subhead{2.1. Reproducing kernels for Riemann surfaces}\endsubhead First, we will  consider
some results of Kazhdan (unpublished), Rhodes \cite {R},  and Jorgenson
and Kramer \cite {JK}.

Let 
$C := \Gamma {\backslash} \bold H$ denote a compact Riemann surface of genus $g(C)\geq 2$, where
$\bold H$ is the upper half-plane and $\Gamma$ a discrete subgroup of $PSL_2(\bold R) $.
 Let $\{C_\alpha\}$ be a cofinal
system of finite Galois coverings of
$C$. Set 
$\Gamma_\alpha:=\pi_1(C_\alpha) .$ Let
$\ell$ be a sufficiently large even integer.
 One can estimate the following function on   $\bold H$:
$$
\eusb F_{{\ell/2}, \alpha}(z, \bar z):= \root {\ell/ 2} \of {(y^{\ell/2})^2}\cdot 
\Bigl(\sum_{j=1}^{n_\alpha}|f_j|^2\Bigl)^{2/\ell} = y^2 \cdot \Bigl(
\sum_{j=1}^{n_\alpha}|f_j|^2\Bigl)^{2/\ell},\tag1 
$$
where $\{f_1, \dots, f_{n_\alpha}\}\subset  S_{\ell }(\Gamma_\alpha)$ is
an orthonormal basis of cuspidal forms  and $y=\roman{Im}(z)$.
For
$\beta =\alpha+1,$ the basis of
$S_{\ell }(\Gamma_\beta)$ is formed by adding forms to the basis of $S_{\ell
}(\Gamma_\alpha).$ In \cite {JK}, the authors considered only cuspidal forms of weight 1 (in our
notation), denoted by
$S_{2 }(\Gamma_\alpha),$ but the argument works for $S_{\ell}(\Gamma_\alpha)$ if we replace the
heat kernel 
$K_\alpha^{(1)}(t,z,w)$ by
$K_\alpha^{({\ell/ 2} )}(t,z,w)$.   
The upshot is that  
$$
\qquad \qquad B_{\bold H, \ell}(z,w) := \lim_\alpha \sum_{j=1}^{n_\alpha}f_j(z)\overline {f_j(
w)}\qquad
\qquad (z, w \in \bold H)
$$
is a holomorphic function in $z$ and antiholomorphic function in $w$.  

\subhead{2.2. Reproducing kernel for $\bold{\tilde X}$}\endsubhead We assume that $X$  satisfies 
the assumptions of Section 1.1.  Let
$X_\gamma$ be a cofinal system of finite Galois coverings of $X.$ Set $U:=\tilde X.$ We may
consider a  natural $\pi_1(X)$-invariant Hermitian metric
$h(.,.)$  on
$\eusb K_U^{\otimes m}$, as in \cite{K, Chap.\, 5.13, 7.1, 5.12}.  By \cite{K, Chap.\, 5.13}, the
natural  Hermitian metric on $K_U$ is equivalent to any other $\pi_1(X)$-invariant Hermitian
metric on $\eusb K_U.$ 
  In Section 1.4, we may
replace complex-valued functions by sections of  relevant bundles.

We consider the bundles 
$
\Cal E_{\gamma, m}:=\eusb K^{\otimes m}
_{X_\gamma}
$
and
$
\Cal E_{U, m}:=\eusb K^{\otimes m}
_U
$ for an appropriate fixed $m$ \/  so that each $\Cal E_{\gamma,m}$ is very ample; such an
$m$ exists  
\cite{K, Chap. 16.5}.
They are equipped with the  
natural Hermitian metrics $h(.,.)$ and
$\|f(z)\|= \sqrt {h(f(z),f(z))},$ as in \cite{K, Chap.\, 5.13, 7.1, 5.12}.          
For each  $\gamma,$ we consider a vector space of integrable forms of weight $m$ on
$X_\gamma$ with a norm that behave well in the tower of coverings:
$$
{1\over{ vol(X_\gamma)}}\int_{X_\gamma} \|f\|^2\rho^{-m} dv.
$$
The union of these vector spaces is a pre-Hilbert  $p\!H_{X,m}$ space  whose
completion is denoted by  $H_{X,m}$. 

\remark{Remark {\rm 2.2.1}}Clearly, each $X_\gamma$ satisfies the above assumptions. Given
 a finite number of points on $\tilde X$, we can assume they lie on a connected
open Riemann surface $Y\subset \tilde X$ which is a Galois covering of a curvilinear section of
$X_\gamma$
 ($\gamma \gg 0$, $\gamma$ depends on the points) by Deligne's theorem \cite{K, Theorem 2.14.1}.
\endremark

We will show that the corresponding
$H_{X,m}$ has a reproducing kernel provided $\pi_1(X)$ is very large (see Definition 2.2.2 below).
We consider
$U\times U$ and write a point of 
$U\times U$ as a pair
$(z, w)$.  The reproducing kernel will be a section of a vector bundle $ p^*_1 \Cal E_{U,m}  \otimes
p^*_2
 \bar \Cal E_{U,m}$ where $p_1$ and $p_2$ are the coordinate projections of $ U\times U.$ Let
$\{u_k\}
\subset p\!H_{X,m}$ be an orthonormal basis constructed similarly to the basis of $\bigcup_\alpha
S_{\ell}(\Gamma_\alpha)$. We set
$$
B_{U,m}(z,w):=\sum p^*_1 u_k(z)\otimes p^*_2 \bar u_k(w)
$$
and verify that the series converges and the sum is a section holomorphic in the first variable.
We take a general hyperplane section $\tau: D \hookrightarrow X$ through a point $x\in X$ 
and consider a standard exact sequences with
$D=m K_X$ ($m$-th power of the canonical divisor):
$$
   0\longrightarrow
\eusb K_X\longrightarrow \eusb K_X(D)
\overset{P.R.}\to
\longrightarrow
\eusb K_D\longrightarrow 0
$$
where $P. R.$ is the Poincar\'e residue map, and $\tau^*(\eusb K_X(D))= \eusb K_D.$    We also
have a natural inclusion
$$
 \Cal O_X(mK)=\Cal O_X(D)\hookrightarrow \eusb K_X(D)= \Cal O_X((m+1)K).
$$
The natural
Hermitian metric on $\eusb K_X(D) $ will induce a Hermitian metric on $\eusb K_D$ as well as on
$\Cal O_X(D)$.  Furthermore, all metrics behave well in the tower of Galois coverings of $X.$ By
trivial induction, we are reduced to the one-dimensional case.  Let $Y$ denote  the inverse image
on
$\tilde X$ of a general curvilinear section
$C\subset X$ through a point of $X.$   It is connected.  One can estimate on $Y$ the function similar
to the expression (1).  It follows that the
series converges pointwise hence uniformly by Dini's theorem, and $ B_{U,m}(z,w)$ satisfies
the condition 1.4(i).
\definition{Definition 2.2.2} With the above assumptions and notation, we say  $\pi_1(X)$ is {\it
very large}\/  if $ B_{U,m}(z,w)$ satisfies the conditions 1.4(i) and 1.4(ii) (sheaf-theoretic version),
for all $m\gg0$.
\enddefinition

Let $Q\in \tilde X$ be an arbitrary point. Consider a coordinate neighborhood of $Q$ in $\tilde X.$
In that neighborhood, the elements of $H_{X,m}$ become functions  and we may consider the
evaluation map $e_Q.$   With 
assumptions of the theorem, $e_Q$  is a continuous linear functional.  So, for   $m \gg 0,$ 
$H_{X,m}$ is a {\it unique  Hilbert space}\/ with reproducing kernel $ B_{U,m}(z,w)$  (see Faraut's
article in
\cite {F-K, pp.\, 5-12}).  We
obtain a holomorphic {\it embedding} \/ of $U$ into the corresponding infinite-dimensional 
projective space by the Schwarz inequality (see Kor\'anyi's article in \cite{FK, p.\,191}). Indeed,
each finite covering is embedded in the corresponding finite-dimensional  projective space.

Thus, we get a metric on $\tilde X$ and a natural isometric and holomorphic embedding
$$
\Upsilon : \tilde X \hookrightarrow \bold {P}(H^*_{X,m}).
$$
\proclaim{Proposition 2.3} With the above notation  and  assumptions of the theorem, the
functional element of the diastasis  generates a strictly plurisubharmonic function on $\tilde X,$
and 
$\tilde X$ is a Stein manifold.
\endproclaim
\demo{Proof} We have two expressions for the diastasis on $\tilde X$ \cite{C, (5)
and (29)}:
$$
 D(\bold p,\bold q)=\Phi(z(\bold p),  \overline {z(\bold p)}) + \Phi(z(\bold q), 
\overline {z(\bold q)}) -\Phi(z(\bold p),  \overline
{z(\bold q)})-\Phi(z(\bold q),  \overline {z(\bold p)}),                   
\tag2
$$
$$
 D(\bold p,\bold q)= \log {\big(\sum^\infty_{\sigma= 0}|\xi^\sigma(\bold
p)|^2\big)\, \big(\sum^\infty_{\sigma= 0}|\xi^\sigma(\bold
q)|^2\big) \over  |\sum^\infty_{\sigma= 0}\xi^\sigma(\bold p)
\,\overline{\xi^{\sigma}( {\bold q})}|^2},                             
\tag3
$$
where $\Phi(z, \overline z)$ is a real-valued analytic potential of our metric \cite{C, Chap.\, 2}
and  
$[\xi^0:
\xi ^1:
\dots]$ are homogeneous coordinates in the projective space $\bold {P}(H^*_{X,m})$. While (3) is a
global expression for the diastasis, (2) is the functional element of the diastasis.   

For a  fixed $\bold p$, from (3), we see that the diastasis on $\tilde X$ is single-valued,
non-negative, and real analytic in
$\bold q$, except at the intersection of  $\tilde X$ in $\bold {P}(H^*_{X,m})$ with the antipolar
hyperplane of
$\bold p$ where $D(\bold p, \bold q)$ becomes infinite, i.e., the denominator in (3)
vanishes. This intersection is clearly either empty or a codimension one subvariety
in $\tilde X$ \cite{C, pp.\,19-20, Theorems 11,\, 12 and Corollaries 1,\, 3}. 

In (2), we observe that the functional element of the diastasis generates a single-valued 
function on
$\tilde X$. Indeed, it generates a function on $Y,$ the preimage on  $\tilde X$ of 
 a  general curvilinear section through a finite number of points of $X_\gamma$ ($\gamma \gg 0$,
$\gamma$ depepends on the chosen points).   The metric and diastasis on $Y$ are induced from
$\tilde X$ \cite{C, Chap.\, 2}.  Here we use that  every  holomorphic bundle  on a noncompact
Riemann surface is trivial, and  we can assume that $Y$ is noncompact and connected by Deligne's
theorem \cite{K, Theorem 2.14.1} since $\pi_1(X)$ is large. The observation, now, follows since
$\tilde X$ is simply connected.

Thus, the functional element of the diastasis  generates a strictly plurisubharmonic {\it function}
on
$\tilde X$. Since
$X$ is compact, the proposition, now, follows from the Oka-Grauert-Narasimhan solution of the 
Levi problem.
\enddemo

\head 3. Comparison of Hilbert spaces\endhead

The aim of this section is to construct a finite-dimensional linear system on $Z'_{i,\gamma}$ (see
3.1)  that together with a pullback of a linear system on $\tilde X$ will give an embedding of
$Z'_{i,\gamma}$  into projective space. Furthermore, such linear systems will be closely related
for different $i$'s $(1\leq i \leq u);$ see Proposition 3.9 below.

\subhead 3.1. Notation and assumptions\endsubhead We keep the notation of Section 1 and
assumptions of the theorem. Furthermore, we assume
 $\phi_{i,\gamma}$ can not be extended to an etale covering of $X.$   We get the following
diagram of base extensions:
$$
\CD
\qquad Z'_{i,\gamma} @>{\phi'_{i,\gamma}}>>  Z_ {i,\tilde X} @>>>  \tilde X \\
\qquad@VV{\xi_{i,\gamma}}V          @VV{\phi_i}V  @VV{\beta}V  \\
  \qquad Z_{i.\gamma}@>\phi_{i,\gamma}>>   Z_i@>>>  X.\\
\endCD
$$

 Let $s\in \bold N$ denote a sufficiently large integer. We consider an arbitrary square
integrable, or $L^2$-integrable, section that is a restriction of a section over
$M_{i,\gamma}$: 
$$
\qquad\omega_{i,\gamma, s}\in H^0(Z_{i,\gamma},\eusb K^{\otimes
s}_{M_{i,\gamma} }),
\qquad \;\;\ 
\|\omega_{i,\gamma,s}\|^2 =  {1\over{ vol(Z_{i,\gamma})}}\int_{Z_{i,\gamma}} 
|\omega_{i,\gamma,s}|^2 {\rho}^{-s}dv <\infty
.$$  
Here $\omega_{i,\gamma, s}$ is viewed as an automorphic form of weight $s$ on $ \tilde Z_i,$ and
$\rho$ arises from $\rho_X.$  The pullback gives a section
$\omega_{i,\gamma, s}'\in H^0(Z'_{i,\gamma}, \xi^*_{i,\gamma}(\eusb K^{\otimes
s}_{M_{i,\gamma}}|Z_{i,\gamma})).
$

We denote by $\Cal A_{\tilde Z_i}$ the space of holomorphic functions on $\tilde Z_i$, and by 
$
B
$ 
the corresponding reproducing kernel.   To simplify notation, we specify only the
first rows of the matrices in Lemma 3.4 below.

\subhead 3.2. Compact exhaustions\endsubhead Let $K_{i,1}\subset \cdots \subset K_{i,t}
\subset\cdots
\subset Z_{i,\tilde X} = Z_i\times_X\tilde X$ be a compact exhaustion of $Z_{i,\tilde X}$ such that
$$
\tilde K_1\subset \cdots \subset \tilde K_t \subset\cdots \subset \tilde X \qquad (\tilde K_t =
\bigcup_{i=1}^uK_{i,t})
$$
is a compact exhaustion of $\tilde X.$ We assume each compact is a union of compact
neighborhoods. Let $\Phi_i\subset \tilde Z_i$ denote the fundamental
region of $Z'_{i,\gamma},$ and
$K_{i,t,\Phi_i}
\subset
\Phi_i$  denote the inverse image of $K_{i,t}.$ For a given $t,$ we can find a sufficiently high finite 
Galois covering
$\beta_\mu :X^\mu \rightarrow X$ such that
$\tilde K_t$ projects one-to-one into $X^\mu$ while $K_{i,t}$ projects into  $Z_{i,X^\mu}
:=\beta^{-1}_\mu (Z_i).$ Similarly, one defines $\Phi^\mu_i$ and $K_{i,t,\Phi^\mu_i}
\subset
\Phi^\mu_i.$
Set
$$
Z^\mu_{i,\gamma}:= Z_{i,\gamma}\times_{Z_i} Z_{i, X^\mu}, \qquad \Sigma^\mu := Gal (\tilde
Z_i/Z^\mu_{i,\gamma})\;\;{\roman {and}}\; \;\Sigma' := Gal (\tilde
Z_i/Z'_{i,\gamma}).
$$

Given an index $i\; (1\leq i\leq u),$  an arbitrary $\epsilon>0,$ and   
$K_{i, t}$
as above,  we can
$\epsilon$-approximate  the Poincar\'e series $ P_{\omega'_{i,\gamma,s}}$ on $ K_{i,t}$
 by a global section of
$\eusb K^{\otimes s}_{\tilde X}$
since $\tilde X$ is a Stein manifold; indeed, $ P_{\omega'_{i,\gamma,s}}$ may be viewed as a
section of $H^0(K_{i,t}, \eusb K^{\otimes s}_{\tilde X} ).$

\subhead 3.3. System of equations\endsubhead Let $x\in\tilde Z_i \subset \bold C^n$ be an arbitrary
point.  Let 
$$
\Gamma : =Gal(Z^\mu_{i,\gamma} /Z_{i,X^\mu})=\{g_1\Sigma^\mu, g_2\Sigma^\mu, \dots,
g_q\Sigma^\mu\}
$$ 
denote
the corresponding Galois group with unit element
$g_1$ in $ Gal(\tilde Z_i/Z_i)$. 

For a {\it fixed}\/ index $i,$ let $J_{\lambda}(x)$ denote the Jacobian of an automorphism  $\lambda$
on
$\tilde Z_i.$ Let
$m= m(\gamma)\in
\bold N$ denote a sufficiently large integer. 

We consider the following system of $q$ equations in $q$ unknowns $Y_1, \dots,Y_q$ at
$x$:
$$
\cases
\qquad J^m_{g_1}(x)Y_1  +\qquad J^m_{g_2}(x)Y_2 +\cdots+\qquad J^m_{g_q}(x)Y_q\, = G_1(x) \\ 
\;\;\;\,J^{m+1}_{g_1}(x)Y_1 +\quad J^{m+1}_{g_2}(x)Y_2 +\cdots +\quad J^{m+1}_{g_q}(x)Y_q\, =
G_2(x)\\
\;\;\; \cdots\cdots\cdots\cdots\cdots\cdots\cdots
\cdots\cdots\cdots\cdots\cdots\cdots\cdots\cdots\cdots\cdots\cdots\cdots \\ 
J^{m+q-1}_{g_1}(x)Y_1 +J^{m+q-1}_{g_2}(x)Y_2
+\cdots+J^{m+q-1}_{g_q}(x)Y_q\,= G_q(x)  \endcases \tag4
$$
where each $ G_v(x)$ $(1\leq v \leq q) $ is an arbitrary  $Gal(\tilde Z_i/Z_{i, X^\mu}) $-automorphic
form of weight $m$ on  $\tilde Z_i,$
i.e., $J^m_\alpha(x)G_v(\alpha x)= G_v(x)$ for $\alpha \in Gal(\tilde Z_i/Z_{i, X^\mu}).$

\proclaim{Lemma 3.4 (lifting)} For almost every
choice  of representatives $g_1, \dots, g_q$ of the
group\/ $\Gamma$
{\rm(}possibly with a finite number of exceptions{\rm)}  and every  $ x\in K_{i,t,{\Phi^\mu_i}},$
\roster
\item"{(i)}" the system {\rm(4)} has a solution obtained by Cram\'er's rule,\/ and 
\item"{(ii)}" suitable solutions $Y_v$'s are functions of \/ $x\in K_{i,t,{\Phi^\mu_i}}$ and $Y_v(x)
=Y_1(g_vx)$ for $v\geq 2,$ where $\{Y_v(.)\}$ denotes a solution of\/ $(4)$ at the corresponding
point.
\endroster
\endproclaim
\demo{Proof}The $D(x):=Det (J^m_{g_1}(x), \dots,  J^m_{g_q}(x)
)$ is a Vandermonde determinant hence 
$$
D(x)= J^{m}_{g_1} (x)\cdots J^{m}_{g_q} (x)\cdot \prod_{j>k}
(J^m_{g_j} (x) - J^m_{g_k}(x)) \quad (J_{g_v g^{-1}_kg_k}(x) \equiv J_{g_vg^{-1}_k}(g_k
x)J_{g_k}(x) ). 
$$
Clearly
$\Sigma^\mu$ is infinite. Given the unit 
$g_1$ and the point $x$ as above, we can choose, in infinitely many ways,  new representatives $g_v
$'s in the corresponding cosets such that
$ J_{g_vg^{-1}_1}(g_1 x)
\neq 1$ \;($2\leq v\leq q$ and $\forall x\in K_{i,t,{\Phi^\mu_i}}$).  The first assertion follows from the
uniform convergence of Poincar\'e series by trivial induction.

To prove (ii), we will show   that  $Y_k(x)=Y_1(g_kx)$ for $2\leq k \leq q$ and suitable solutions. We
consider a linear system like (4)  with $x$ replaced by $g_kx$ in (4) and the determinant $D(g_kx) .$ From
the identity
$J_{g_j g_k}(x) \equiv J_{g_j}(g_kx)J_{g_k}(x)$, we obtain 
$$
D(g_kx) = Det(J^m_{g_1g_k}(x)J^m_{g_k}(x)^{-1}, \dots,  J^m_{g_qg_k}(x)J^m_{g_k}(x)^{-1}) 
$$
with $J_{g_1g_k}(x)J_{g_k}(x)^{-1}\equiv1.$ Clearly $D(g_kx)$ is the determinant of the
system with the right-hand side $\{G_v(g_kx)\}$ ($1\leq v\leq q$). Hence
$$
D'(x): = Det(J^m_{g_1g_k}(x), \dots,  J^m_{g_qg_k}(x))
$$
is the determinant of a
system   with the right-hand side
$\{G_v(g_kx)J^m_{g_k}(x)\}$.
 But $D'(x)$ is also the determinant of the system, similar
to  (4), in the $q$ unknowns
$Y_k,\dots, Y_1,\dots,$ with right-hand side
$\{G_v(x)\}$. Clearly $G_v(x)= G_v(g_kx)J^m_{g_k}(x)$ for  $1\leq v \leq  q $. It follows that
$Y_k(x)=Y_1(g_kx)$ for suitable solutions. 
\enddemo

\proclaim{Lemma 3.5} We fix\/ $i$ and the compact $K_{i,t}\subset X^\mu$
$(see\; 3.2)$. With    holomorphic on $\tilde Z_i$ automorphic forms $G_1(x), \dots, G_q(x)$
$(see\; 3.3)$,       let
 $F(x)$ denote  a suitable solution of\/
$(4)$ on the compact $K:= K_{i,t,{\Phi^\mu_i}}$ obtained in Lemma {\rm 3.4.}
Set $F(x) = 0$ on $\tilde Z_i\backslash  K.$ Let $L_m^2(\tilde Z_i,\Sigma^\mu, dv) :=L_m^2(\tilde
Z_i,\Sigma^\mu, \rho_X, dv)$ denote the standard Hilbert space (with reproducing kernel).
 Then 
\roster
\item"{(i)}" $P_{F,K}(x):= \sum_{\sigma\in \Sigma^\mu} F(\sigma x)J_\sigma^m(x)\in
L_m^2(\tilde Z_i,\Sigma^\mu, dv) ;$ 
\item"{(ii)}" we obtain a holomorphic on
$\tilde Z_i$ automorphic form of weight $m$ with respect to $\Sigma^\mu$  by applying  the
Bergman projection 
$$
\qquad \quad (\beta_mP_{F,K})(x):= \int_{\Phi^\mu_i} B_m(x,\xi )P_{F,K}(\xi)\rho^{-m}(\xi) dv_\xi
\,\in L_m^2(\tilde Z_i,\Sigma^\mu, dv)\bigcap \Cal A_{\tilde Z_i};
$$
\item"{(iii)}"$P_{F,K}(x)=F(x)$ on $K,$ and $(\beta_mP_{F,K})(x)$ is uniquely determined.
 \endroster
\endproclaim
\demo{Proof} Since $F(x)$ is bounded, the series in (i) converges absolutely and uniformly, and
$P_{F,K}(x)$ is a measurable form in $L_m^2(\tilde Z_i,\Sigma^\mu,
dv).$ Finally, the uniqueness follows  from  Cram\'er's rule and the well-known definition of relative
Poincar\'e series. Recall that the latter definition is independent of the choice of
representatives.
\enddemo 

\subhead 3.6\endsubhead For a fixed \/ $i$ $(1\leq i\leq u) ,$ we consider a finite basis of extendable
sections  
$$
\Lambda_{i,1},\dots, \Lambda_{i,\ell} \in
H^0(Z_{i,\gamma}, \eusb K^{\otimes m}_{M_{i,\gamma}})
$$
that defines an embedding of $ M_{i,\gamma}$ and $ Z_{i,\gamma}$ into projective space. We can
assume that $m=m(\gamma) \gg 0$ and those sections are $L^2$-integrable.  
By Lemmas 3.4 and 3.5, we can find a linearly independent subset of
$$
H^0(Z'_{i,\gamma}, \xi^*_{i,\gamma}(\eusb K^{\otimes
m}_{M_{i,\gamma}}|Z_{i,\gamma}))
$$
that defines an embedding into projective space, consisting of pullbacks of the orthonormal
basis  of the Hilbert space $ H_{ X,m}$ supplemented by the pullback of the maximal (finite)
number of  sections
$$
\Omega_{i,1}, \dots,\Omega_{i,p}, \dots\,\in \{\Lambda_{i,1},\dots, \Lambda_{i,\ell}\}
$$ 
that are not
arising from $ H_{ X,m}$, and, if possible, a finite number of {\it new}\/ $L^2$-integrable
sections.
 The latter sections are coming from $Z_j$'s for
$j\neq i$; see (3.7) below. 
For each $i$, those  sections which are {\it not}\/ arising from
$ H_{ X,m}$ generate a finite-dimensional Hilbert space, denoted by $V_{i,\gamma}$, and
${\roman {dim}}\, V_{i,\gamma}$   is independent of
$i$. For $1\leq k\leq u $,  let $ L_m(\tilde Z_k,\Sigma', dv)$ denote the completion of  $ \bigcup_\mu
L_m^2(\tilde Z_k,\Sigma^\mu, dv) .$

We consider a  sequence $\{\epsilon_t\}$ of small positive real numbers with $\lim
\epsilon_t=0$. Let 
$$
\CD
\Cal Z_i= M_{i,\gamma}\times_X \tilde X @<<<   M_{i,\gamma}\times_X \tilde
X\times_X M_{j,\gamma} @>>>\quad\quad M_{j,\gamma}\times_X \tilde X=\Cal Z_j\\
\quad\quad \;\;@VVV @VVV\!@VVV  \\
 \Cal Z^t_i =  M_{i,\gamma}\times_X   X^t@<f^t_i<<   M_{i,\gamma}\times_X
 X^t\times_X  M_{j,\gamma}@>f^t_j>>\quad\quad M_{j,\gamma}\times_X  X^t=\Cal Z^t_j\\
\endCD
$$
be a natural  diagram  , where $1\leq i,j\leq u$   and $X^t
\rightarrow X$ is a suitable finite Galois covering, i.e., the compact $\tilde K_t$ maps one-to-one
onto a compact $K_t\subset X^t.$ The maps $f^t_i$ and $f^t_j$ are finite and generically etale.

Set 
$\Omega_i:= \Omega_{i,1}.$ Precisely, we have the following correspondence.

\subhead {3.7. Approximation}\endsubhead Given the  element  $\Omega_i$,  now   we will
produce an element $\Omega_j \in L_m(\tilde Z_j,\Sigma', dv) $ for $1\leq j \leq u$. By abuse
of notation, we denote  by
$\Omega_i$ the  pullbacks of
$\Omega_i$ to $\Cal Z_i$ and $\Cal Z^t_i$.
Our compact exhaustions (see (3.2)) give  compact exhaustions of $\Cal Z_i$.
For a suitable $t$, the corresponding compact set  descents to $\Cal Z^t_i$. With the above
notation, let
$P_{\Omega_i}(x)$ denote the corresponding Poincar\'e series (see\; 3.2). The Poincar\'e series can be
$\epsilon_t$-approximated on
$K_{i,t}$ by global sections of\/ $\eusb K^{\otimes m}_{\tilde X}$. We also consider
$P_{\Omega_ih^w_i}(x)  $, where 
$h_i$ is the linear form defining $X\backslash Z_i$. Let 
$$
G_{t,1}h_i^0\in H^0(\tilde X, \eusb K^{\otimes m}_{\tilde X}),\dots, G_{t,q} h_i^{q-1}\in H^0(\tilde X,
\eusb K^{\otimes {m+q-1}}_{\tilde X})
$$
be the corresponding $\epsilon_t$-approximating sections on the compact $K_{i.t}$, where each 
$G_{t,v}\in H^0(\tilde X, \eusb K^{\otimes m}_{\tilde X})$.
 By   lemmas, we obtain automorphic
forms 
$$
\omega_{j, t}\in L_m^2(\tilde Z_j,\Sigma^t, dv)\bigcap \Cal A_{\tilde
Z_j}, \qquad \lim_{t \to
\infty} \omega_{j, t}=\Omega_j \in  \Cal A_{\tilde Z_j},
$$
where the convergence is in the Hilbert space $ L_m(\tilde Z_j,\Sigma', dv) .$ The convergence
follows from the assumptions on $\Omega_i$, the uniqueness (see  (3.4)-(3.5)), and the
diagram in (3.6) that allows  to compare $\Cal Z^t_i$ and $\Cal Z^t_j$.   Thus
$\Omega_j$ is the desired element.

\subhead{3.8. Notation}\endsubhead Let $H_{i,\gamma,m}:=V_{i,\gamma} \oplus H_{X,m}$ denote a
Hilbert space, and
$A_{\gamma, m}$ the corresponding infinite-dimensional complex linear space (we have applied a
forgetful functor). Let
$\bold P(A_{\gamma, m}^*)$ denote the corresponding projective space. We will consider only the
linear transformations of
$\bold P(A_{\gamma, m}^*)$ that leave $ \bold P(H_{X,m}^*) $ fixed.

We have proved the following
\proclaim{Proposition 3.9}We keep the above notation and assumptions of the theorem, and fix 
$\gamma$ and
$m=m(\gamma)\gg 0.$ For all
\/
$i$'s
$(1\leq i\leq u)$, there is a natural one-to-one linear correspondence between
  the Hilbert spaces \ 
$H_{i,\gamma,m}$  that induces an isomorphism on the Hilbert subspace $H_{X,m}\subset
H_{i,\gamma,m }.$ Moreover,
${\roman {dim}}\, V_{i,\gamma}$ is finite and independent of\/ $i.$
\endproclaim

\head{4. Monodromy}\endhead
\subhead{4.1. Notation and Assumptions}\endsubhead We keep the notation of Section 1 and
assumptions of the theorem.   Let $B_i(\tilde X)$ denote  the inverse image of the hyperplane
$X\backslash Z_i$ on
$\tilde X$. Each
$\tilde X{\backslash}B_i(\tilde X)$ is connected and the natural map 
$\pi_1 (Z_i) \rightarrow \pi_1 (X)$ is an epimorphism. We get a natural holomorphic
surjection 
$\phi'_{i,\gamma}: Z'_{i,\gamma} \rightarrow Z_{i,\tilde X}.$ Set $U_i:=Z'_{i,\tilde X}$
 and $g_i:=\phi'_{i,\gamma}.$

We fix an arbitrary point $\tilde x \in \tilde X \backslash (\bigcup_iB_i(\tilde X))$. Let $f(0)$
be a unique natural germ of holomorphic map from $\tilde X$ into $ U_1 \subset
\bold P(A^*_{\gamma,m})$ sending $\tilde x$ to a point $\tilde u \in U_1$ such that
$g_1(\tilde u)=\tilde x$. Finally, we fix the inclusion $\tilde Z_1 \subset \bold C^n$.

\proclaim{Proposition 4.2} With  the above notation and assumptions of the theorem, there exists a
unique holomorphic inclusion
$$
f_1: \tilde X \hookrightarrow \tilde Z_1\subset \bold C ^n.
$$ 
In particular, $\tilde X$ is holomorphically equivalent to a bounded domain in $\bold C^n.$
\endproclaim

\demo{Proof} We proceed in two steps. First, we will construct a holomorphic 
map
$ \Upsilon_\gamma:\tilde X\rightarrow \bold P(A^*_{\gamma, m})$ and a 
holomorphic inclusion $\tilde X\backslash B_1(\tilde X) \hookrightarrow\tilde Z_1 \subset \bold
C^n.$
 Then we will map $\tilde X$ into $\bold C^n$ and verify that we get an inclusion.

{\it {Step $1$}}.  We take an arbitrary point $P_1\in \tilde X$ and connect it to $\tilde x$ by a
continuous path $P_t\subset \tilde X \ (0\leq t\leq1, P_0=\tilde x)$. We
will show that there exists one and only one family  $f(t) \ (0\leq t\leq1)$ of
germs of holomorphic maps from $\tilde X$ into $\bold P(A^*_{\gamma, m})$ sending $P_t$ to
$f(t)(P_t)$ such that $f(t)$ induces $f(t^*)$, when $P_{t^*}$ is near enough to $P_t$, and
$f(0)$ is our fixed germ; moreover,  $f(t)(P_t) \in U_1$ and
$g_1\!\cdot\! f(t)(P_t)=P_t  $ provided
$P_t$ lies under a point of $U_1$.   Such a family is called compatible.

We claim 
that if $f(t) \ (0\leq t<1)$ is a compatible family of holomorphic maps which is
open on the right side then we can complete it by $f(1)$ to a compatible family $f(t) \ (0\leq
t\leq1)$.

Let $\phi =\phi^{P_1}$ be a holomorphic map of a small neighborhood $\Cal V\subset
\tilde X$ of $P_1$ into
$\bold P(A^*_{\gamma, m})$ corresponding to a germ sending $P_1$ to a point in $U_i$ such
that
$\phi (\Cal V) \subset U_i$ and $g_i  \cdot\phi (\Cal V) = \Cal V$. Such a germ always
exists.

First, we will assume that $i=1$, i.e., 
$P_1\in \tilde X\backslash B_1(\tilde X)$ . Hence $P_1$ lies under a point $Q\in U_1$ and
$\phi (\Cal V)$ is a small connected neighborhood in $U_1$ that projects onto
$\Cal V$, i.e., $g_1\cdot\phi$  is an identity on $\Cal V$.  

Let $P_s$ be a point of the path which is near enough to
$P_1$ such that 
$P_t \in\Cal V$ for $s\leq t\leq1$.   Let $\Cal W\subset \Cal V$
be a neighborhood of $P_s$ such that
$f(s)$ maps it into $\bold P(A^*_{\gamma, m})$. In this case, $f(s)\cdot \phi^{-1}$ gives a
holomorphic map from $\phi(\Cal W)$ onto $f(s)(\Cal W)$.

The key point is that the map $f(s)\cdot \phi^{-1}$ can be {\it extended uniquely}\/ to a
global holomorphic linear map \/
$\sigma$ \/ of $\bold P(A^*_{\gamma, m})$   because $\phi (\Cal W)$ spans linearly 
 $\bold P(A^*_{\gamma, m})$.  

Moreover, $\sigma$ acts as a deck transformation of
$U_1$ over  $\tilde
X\backslash B_1(\tilde X)$   because it is a deck transformation over a small neighborhood in $\tilde
X\backslash B_1(\tilde X)$.  We get
$\sigma (U_1) =U_1$ and 
$g_1\cdot \sigma (a) =g_1(a)$ for every $a \in U_1$.
It is a simple matter to see that the germ
$f(1)$ determined by $\sigma\cdot \phi$ at  $P_1$  satisfies our demand.

Now, we will assume that $P_1\in B_1(\tilde X)$ and $P_1$ lies under a point $Q\in 
U_i,$ where $i\neq1$. A similar argument applies to obtain a local map as well as a
global holomorphic linear map $\sigma$.

 This establishes our claim, i.e., we can complete the family on the right.
Because $\tilde X$ is simply connected we get a holomorphic map  
$\Upsilon_\gamma:\tilde X
\rightarrow
\bold P(A^*_{\gamma,m})$.  Since $\pi_1(Z_i)$ is residually finite, we get a
holomorphic inclusion 
$$
\tilde X\backslash B_1(\tilde X) \hookrightarrow\tilde Z_1 \subset \bold C^n.
$$

{\it Step $2$}\, (inclusion into $\bold C^n$). Because
$B_1(\tilde X)\subset \tilde X$ is a holomorphic subset, the above inclusion
 has a unique
holomorphic extension
$
f_1: \tilde X \rightarrow\bold C^n. 
$

We will show that $f_1$ is a  holomorphic inclusion hence $\tilde X$ is a
bounded domain.
If dim $\tilde X =1$, i.e., $n=1$ then $f_1$ is an open bijection onto its image by
the structure theorem for holomorphic functions in one variable. We will reduce the
general case to the one-dimensional case. 

We suppose that $f_1$ is {\it not}\/ a one to one map onto its image and derive a
contradiction. Let   $q\in f_1(\tilde X)\subset \bold C^n$ be an exceptional point, i.e.,
$f_1^{-1}(q)$ contains more than one point.                
By Osgood's theorem,  there is a point $Q\in  f_1^{-1}(q)$ which is {\it not}\/ an isolated point
in $f_1^{-1}(q)$ unless an open
neighborhood of $q$ in $\bold C^n$ is covered exactly $k$ times by an open neighborhood
of $Q$ in $\tilde X$. In the latter case, we get $k=1$, a contradiction.  In the former case, we
will first assume that dim\,$X=2.$

{\it Case}\,: dim\,$X=2$.  Then $f_1^{-1}(q) = B_1(\tilde X)$ since $B_1(\tilde X)$ is a nonsingular
connected curve. Let $C\subset X$ be a general hyperplane section. Its preimage in $\tilde
X$, denoted by $\beta^{-1}(C)$, will intersect $B_1(\tilde X)$ in an infinite number of
points, and an infinite number of points of $\beta^{-1}(C)$ will go to the point $q$ under the
map $f_1$. 

We consider the following two diagrams: 
$$
\CD
d\in\Delta@<f_{1,C}<<  \tilde C   @. \qquad \qquad  \Delta @<f_{1,C}<<  \tilde C \\
\;\;\;\;\;\;@V{\tau}VV          @VV{\nu}V \qquad \qquad @V{\tau}VV          @VV{\nu}V\\
 \tilde u \in\tilde Z_1  @.  \tilde X  @. \qquad \qquad \bold C^n @<f_1<<\beta^{-1}(C).\\
\;\;\;\;\;\;@V{\alpha}VV @VV{\beta}V   \\ 
c\in Z_1 @>{}>> X @. @.
\endCD
$$  
Here the disk $\Delta \subset \bold C$ denotes the universal covering of $C\cap Z_1$, as well
as of its preimage under $\alpha$ by \cite{K, Theorem 2.14.1}, and $\tau (\Delta)\subset \tilde
Z_1\subset
\bold C^n$. The existence of the holomorphic inclusion
$f_{1, C}
$ was just established. We will show that the map $f_1$, restricted to $\beta^{-1}(C)$, is
also an inclusion {\it into}\/ $\tilde Z_1$. 

We claim that the second diagram is commutative provided we have fixed a point
$c\in C\cap Z_1$ and have chosen a point $\tilde u  \in\tilde Z_1$ over $c$ and a point
$d\in\Delta$ over $\tilde u$ together with their respective small complex neighborhoods.
Indeed, on $\tilde C\backslash B_1(\tilde C)$, we get
$$
\tau \cdot f_{1,C} = f_1 \cdot \nu
$$
by the constructions of the maps $f_1$ and $f_{1, C}$ employing analytic continuations
along paths, as in Step 1; $B_1(\tilde C)$ is defined similarly to $B_1(\tilde X).$ The equality holds
on $\tilde C$ by continuity. It follows that the point $q$ actually belongs to $\tilde Z_1$. 

Furthermore, $\tau$ is clearly an open map from $\Delta$ onto $\tau(\Delta)$ hence $\tau
\cdot f_{1, C}$ is an open map from $\tilde C$ onto $\tau\cdot f_{1, C} (\tilde C)$. Its
image contains a small open disk about the point $q$, say $\Cal V_q \subset \tau\cdot f_{1,
C} (\tilde C)$. We have a holomorphic inclusion
$ 
f_1^{-1} : \Cal V_q\backslash q \hookrightarrow \beta^{-1}(C).
$
We also have a holomorphic map
$$
f_2: \beta^{-1}(C) \longrightarrow \tilde Z_2 \subset \bold C^n
$$
where $f_2$ and $\tilde Z_2$  are similar to $f_1$ and $\tilde Z_1$. Clearly we may
assume that $f_2\cdot f_1^{-1}|\Cal V_q\backslash q$ is an inclusion into a bounded
domain in a copy of $\bold C^n$ since $C$ was a {\it general}\/ hyperplane section. It
follows that the map $f_1^{-1}|\Cal V_q
\backslash q$ can be extended to the point $q$. Hence $f_1^{-1} (q)$
 can not contain more than one point, a contradiction.

{\it Case}\,: dim\,$X\geq 3$. We take a general $2$-dimensional linear section $Y \subset
X$ that intersects the image of $f_1^{-1}(q)$ in $X$ in at least two distinct points; this is
possible since not every point of $f_1^{-1}(q)$ is isolated in $f_1^{-1}(q)$. By the
Lefschetz-type theorem, the natural map $\pi_1(Y) \rightarrow \pi_1(X)$ is an isomorphism.
Furthermore, $Y$ satisfies all the assumptions of the theorem and dim\,$Y=2$, a
contradiction.

This proves the proposition and the theorem.
\enddemo

\remark{Remark {\rm 4.3 (Campana)}}The universal covering of a
sufficiently general ample divisor $X$ in a  simple Abelian threefold  is not a bounded
domain, since it follows  from \cite{LS} that all bounded holomorphic functions on $\tilde X$ 
are  constants; $\pi_1(X)$  is large and residually finite, and $ \eusb K_X$ is ample.
\endremark

\smallskip
{\it Acknowledgment.} I would like to thank J\'anos Koll\'ar for
helpful conversations, and Peter Sarnak for bringing the reference \cite
{R}  to my attention  in 2004.
I am especially grateful to Fr\'ed\'eric Campana for pointing out an error in an earlier version of 
the note.

\enddocument